\documentclass[psamsfonts]{conm-p-l}

\usepackage{amssymb,amsmath,latexsym,amsfonts}
\usepackage[all]{xy}

\newtheorem{Prop}{Proposition}[section]
\newtheorem{Thm}[Prop]{Theorem}
\newtheorem{Lemma}[Prop]{Lemma}

\numberwithin{equation}{section}

\newcommand{\Zz}{\mathbb{Z}}
\newcommand{\Aa}{\mathbb{A}}
\newcommand{\Cc}{\mathbb{C}}
\newcommand{\Gg}{\mathbb{G}}
\newcommand{\Rr}{\mathbb{R}}

\newcommand{\mclus}{\mathrm{clus}}
\newcommand{\mvol}{\mathrm{vol}\,}
\newcommand{\mconv}{\mathrm{conv}\,}
\newcommand{\mExt}{\mathrm{Ext}}
\newcommand{\mHom}{\mathrm{Hom}\,}
\newcommand{\mGL}{\mathrm{GL}}
\newcommand{\mand}{\mbox{ and }}

\newcommand{\cT}{\mathcal{T}}
\newcommand{\cR}{\mathcal{R}}

\newcommand{\udim}{\mathrm{\underline{dim}\,}}

\newcommand{\lra}{\longrightarrow}

\begin{document}

\title[Tilting modules, polytopes, and Catalan numbers]{Tilting Modules over the Path Algebra
of Type $\Aa$, Polytopes, and Catalan Numbers}

\author{Lutz Hille}
\address{Mathematisches Institut, Universit\"at M\"unster, Einsteinstrasse 62, D-48149 M\"unster, Germany}
\email{lutz.hille@uni-muenster.de}

\dedicatory{Dedicated to Helmut Strade on the occasion of his 70th birthday}



\maketitle

\date{}


\begin{abstract}
It is well known that the number of tilting modules over a path algebra of type $\Aa_n$ coincides
with the Catalan number $C_n$. Moreover, the number of support tilting modules of type $\Aa_n$
is the Catalan number $C_{n+1}$. We show that the convex hull $C(\Aa_n)$ of all roots of a root
system of type $\Aa_n$ is a polytope with integral volume $(n+1)C_{n+1}={2n\choose n}$. Moreover,
we associate to the set of tilting modules and to the set of support tilting modules certain polytopes
and show that their volumes coincide with the number of those modules, respectively. Finally, we
show that these polytopes can be defined just using the root system and relate their volumes so
that we can derive the above results in a new way.
\end{abstract}


\section{Introduction}

We consider a quiver $Q$ of type $\Aa$ with $n$ vertices, more details can be found in Section 2.
An indecomposable representation of $Q$ can be identified with an interval $[i,j]$ representing the
support of the dimension vector $(0,...,0,1,\ldots,1,0,...,0)$. We associate to $Q$ two series of
polytopes, the $C$-- and the $P$--series, which both come in three versions. The first series of
polytopes $C^+(Q)\subset C^{\mclus}(Q)\subset C(Q)$ consists just of the convex hulls of
certain roots in the root system of type $\Aa$. Thus these polytopes are defined independently
of the orientation of $Q$.

All modules in this note should be understood as modules over the path algebra of type $\Aa$. A
tilting module is a particular module satisfying certain genericity properties. We always identify
modules over the path algebra with representations of the corresponding quiver.

The second series of polytopes $P^+(Q)\subset P^{\mclus}(Q)\subset P(Q)$ consists of the union
of certain simplices $\sigma_T$ associated to tilting modules $T$ or some generalizations of these.
Note that the first series of polytopes only depends on the underlying graph of $Q$, whereas the
second one depends on a chosen orientation of the quiver. However, we will show that the polytopes
in the second series also only depend on the underlying graph of $Q$. Thus we will write $C(Q)$ or
$C(\Aa_n)$ interchangeably, but we have to write $P(Q)$ until we haven proven that the latter definition
is independent of the orientation of the quiver.

The principal goal of this paper is to compare both types of polytopes. In fact, we show that they
coincide (see Theorem \ref{Thmequal}). Moreover, we obtain the number of certain versions of tilting
modules as the volume of the corresponding polytope, where we use a certain normalization of the
euclidean volume. A second aim of this paper is to further simplify the counting by passing to tilting
sequences. This approach is explained in Section 5. Here we consider an additional order on the
indecomposable direct summands, and then the counting gets even easier, namely, we just obtain
a bijection between the tilting sequences and a certain symmetric group (see Theorem \ref{ThmTiltSequ}).

In fact, in this paper we do not use any representation theory other than some geometric interpretations
of representation-theoretic notions. For any dimension vector $d$ there is exactly one rigid (or generic)
module that is dense in the corresponding representation space (see Section \ref{sectArep} or also
\cite{BaurHille} for more details), and we will work with this rigid module.

Note that we use the integral volume in this paper (see also Section 6 for more explanations), that
is, the volume $\mvol\,\Delta$ of any simplex $\Delta$ generated by an integral basis is $1$. Thus
our volume is just $n!$ times the usual euclidean volume. The advantage of this definition is that the
volume of any lattice polytope is an integer.

To be precise we define certain positive numbers (not depending on the orientation by Theorem
\ref{Thmequal})
$$
t^+(Q)=\mvol P^+(Q)\,,\quad t^{\mclus}(Q)=\mvol P^{\mclus}(Q)\,,\quad t(Q)=\mvol P(Q)\,.
$$

In order to get an interpretation of these numbers, we need to consider several versions of tilting
modules. Note that a tilting module $T=\oplus_{i=1}^nT(i)$ is a direct sum of $n$ pairwise non-isomorphic
indecomposable modules $T(i)$ satisfying $\mExt^1(T,T)=0$. For an indecomposable representation
$[i,j]$ the support is just the interval $[i,j]$, for a direct sum of such modules the support is just the
union of the supports of the indecomposable direct summands. Thus, the support of a module (or a
representation) is the support of its dimension vector (see Section 2 for details). A support (or cluster)
tilting module is a module $T$ that is a tilting module if restricted to its support. A $2$--support tilting
module is a module $T$ together with a decomposition $T=T^+\oplus T^-$ such that both $T^+$ and
$T^-$ are support tilting modules and the supports are a disjoint union of the vertices of the quiver
$Q$ of type $\Aa_n$. Thus, any 2-support tilting module $T$ defines a subset $I$ of the set $Q_0$
of vertices of $Q$ such that $I$ is the support of $T^+$ and $Q_0\setminus I$ is the support of $T^-$.
We define $\cT^+(Q)$, $\cT^{\mclus}(Q)$, and $\cT(Q)$, respectively, as the set (of isomorphism
classes) of tilting modules, the set (of isomorphism classes) of support tilting modules, and the set
(of isomorphism classes) of 2-support tilting modules. By using the results in \cite{VolHille} we obtain
the following interpretation of the volumes of the polytopes associated to tilting modules and their
generalizations.

\begin{Thm}\label{Thmnumbervol}
$$
t^+(Q)=\sharp\cT^+(Q)\,,\quad t^{\mclus}(Q)=\sharp\cT^{\mclus}(Q)\,,\quad t(Q)=\sharp\cT(Q)\,. 
$$
\end{Thm}

Note that the classification of tilting modules $\cT^+(Q)$ over quivers $Q$ of type $\Aa_n$ is well
known. A description using trees for the directed orientation can be found in \cite{VolHille}. This leads
to a recursion formula for $t^+(Q)$. In particular, the recursion formula in Theorem \ref{Thmcompare}
is the same as the recursion formula for the number of $3$--regular trees with $n+1$ leaves and one
root (see Section \ref{sectDyck}).

The key point for the correspondence in type $\Aa$ is Theorem \ref{Thmequal} that does not generalize
to the other Dynkin quivers apart from type $\Cc$ and is certainly wrong for euclidean and wild
quivers. Let us briefly comment on the non-simply laced types. There is the notion of a path algebra,
where we use two fields, one being a field extension of degree two (respectively, three for Dynkin
type $\Gg_2$) of the other one. This is certainly more technical and will be explained in our forthcoming
paper \cite{PolDynkHille} in detail. We also note that most of the combinatorics related to the Catalan
numbers is already known for type $\Aa$, but in all other cases it is unknown.

The next theorem shows that for type $\Aa$ the polytopes in the $P$--series are independent of the
chosen orientation. Thus we can simply write $P(\Aa_n)$ instead of $P(Q)$.

\begin{Thm}\label{Thmequal}
$$
C(Q)=P(Q)\,,\quad C^{\mclus}(Q)=P^{\mclus}(Q)\,,\quad C^{+}(Q)=P^{+}(Q)\,.
$$
\end{Thm}

By using decompositions of the polytopes we get several recursion formulas relating the three different
polytopes. Note that $t(\Aa_1)=t^{\mclus}(\Aa_1)=2t^+(\Aa_1)=2$. Consequently, already these
formulas determine the numbers $t(\Aa_n)$, $t^{\mclus}(\Aa_n)$, and $t^+(\Aa_n)$ uniquely, by
induction. However, we can also compute $t(\Aa_n)={2n\choose n}=(2n)!/n!n!$ directly, which gives
even more ways to determine these numbers.

The next result is  a standard decomposition that holds for any quiver $Q$. In Section \ref{sectArep}
we will present two examples illustrating these formulas.

\begin{Thm}\label{Thmcompare}
$$
t(\Aa_n)=\sum_{I\subseteq Q_0} t^+(\Aa_n|_I) t^+(\Aa_n|_{Q_0\setminus I})\,,
$$
$$
t^{\mclus}(\Aa_n)=\sum_{i\in Q_0}t^+(\Aa_n|_{Q_0\setminus\{i\}})=\sum_{i=0}^{n-1}t^+(\Aa_i)t^+(\Aa_{n-1-i})\,.
$$
\end{Thm}

For the $C$--series of polytopes we can prove the following formulas directly by determining the facets
and their volumes, where the polytope $C(\Aa_n)$ is the simplest one to consider.

\begin{Thm}\label{Thmcompute}
$$
\mvol C(\Aa_n)={2n\choose n}=(n+1)C_n\,,\quad\mvol C^{\mclus}(\Aa_n)={2n+2\choose n+1}/(n+2)=C_{n+1}\,,
$$
$$\mvol C^{+}(\Aa_n)={2n\choose n}/(n+1)=C_n\,.
$$
\end{Thm}

In type $\Aa$ we have one dimension vector $(1,1,\ldots,1)$ corresponding to the interval $[1,n]$.
This provides us with yet another recursion formula which is easy to prove for the directed orientation.

\begin{Thm}\label{Thmcompute1}
$$
t^+(\Aa_n)=\sum_{i=1}^n t^+(\Aa_n|_{Q_0\setminus\{ i\}})\,.
$$
\end{Thm}

Since there are several ways to prove these theorems, we give a short outline of their proofs. One
way to prove the results uses induction over the facets and another way uses the fact that a quiver
of type $\Aa$ has precisely one sincere root. However, by using a purely combinatorial approach,
we obtain a different proof by simply using counting arguments as explained in  Stanley's book on
enumerative combinatorics \cite{Stanley}.

We also like to mention that the same idea works for arbitrary Dynkin quivers. However, in this case
the polytope $P(Q)$ is not convex, and thus does not coincide with $C(Q)$, which is the convex hull
of all the roots of the root system. This makes the above formulas more complicated, and for the
details we refer the reader to our forthcoming paper \cite{PolDynkHille}. We would also like to point
out that the number of tilting modules for arbitrary Dynkin quivers has recently been computed in
\cite{Ringeletal} by completely different methods. Moreover, the polytope $P(Q)$ (and its variations)
can be defined even for arbitrary infinite quivers $Q$. This way we would get a simplicial complex
that is a triangulation of a certain quadric in the real Grothendieck group $K_0(Q)_{\Rr}$ of the
category of representations of $Q$.

The outline of the paper is as follows. After the introduction we start with some basic representation
theory in Section 2. Here we only recall a few facts that are well known. Details and further references
can be found in \cite{BaurHille}. In Section 3 we prove the first three theorems and in Section 4 we
proceed with the last two theorems. In Section 5 we modify the problem slightly. Instead of tilting
modules we will consider tilting (or full strongly exceptional) sequences. Finally, in the last section
we will give some combinatorial interpretation of the results using Stanley's book \cite{Stanley}.
\medskip

\noindent {\bf Acknowledgments.} This work started during a stay of the author in Bielefeld at the
SFB 701 'Spectral Structures and Topological Methods in Mathematics'. He would like to thank Henning
Krause for the invitation and the stimulating working conditions. Moreover, this work was supported
by SPP 1388 'Representation Theory'. The author is also indebted to Friedrich Knop for several hints
concerning the combinatorics of the Catalan numbers and to Claus Michael Ringel for discussing further
combinatorial aspects of the Catalan numbers. Finally, he is grateful to Karin Baur, J\"org Feldvoss,
and the referee for many helpful comments on various drafts of this paper.


\section{Representations of $\Aa_n$}\label{sectArep}

In this section we always consider a quiver $Q$ of type $\Aa_n$, that is, a Dynkin diagram of type
$\Aa_n$, where we choose an orientation of any edge between the vertices $i$ and $ i + 1$. For
every oriented edge $\alpha\in Q_1$, called arrow, we define its starting point to be $s(\alpha)$ and its
terminal point to be $t(\alpha)$. A representation of $Q$ consists of $n$ finite dimensional vector spaces
$V_i$ over a fixed field $k$, where $i=1,\ldots,n$, with linear maps $V(\alpha):V_{s(\alpha)}\lra
V_{t(\alpha)}$. Note that either $s(\alpha)+1=t(\alpha)$ or $s(\alpha)-1=t(\alpha)$. The dimension
vector $\udim V$ of the representation $V=\{V_i\}$ is defined as $\udim V=(\dim V_1,\ldots,\dim V_n)$.
A {\sl sincere root} is a dimension vector with $\dim V_i\not= 0$ for all vertices $i$. The set of all
representations of a quiver with the usual homomorphisms forms an abelian category of global
dimension one that has enough projective and also enough injective representations. Note that
any representation $V$ has a short projective resolution $0\lra P^1\lra P^0\lra V\lra 0$ and
$\mExt^1(V,W)$ is defined as the cokernel of the induced map $\mHom(P^0,W)\lra\mHom(P^1,W)$.
The kernel of this map is the set of all homomorphisms $\mHom(V,W)$. A representation $V$ is
{\sl simple} if it has no proper subrepresentations, it is {\sl indecomposable} if it has no non-trivial
decomposition into a direct sum of two representations, and it is {\sl rigid} if $\mExt^1(V,V)=0$.
The latter condition has also a natural interpretation in the space $\cR(Q;d)$ of all representations
of dimension vector $d$ consisting just of all possible linear maps
$$
\cR(Q;d)=\bigoplus_{\alpha\in Q_1}\mHom(k^{d(s(\alpha))},k^{d(t(\alpha))})\,.
$$
The group $G(d)=\prod\mGL(d_i)$ acts on $\cR(Q;d)$ via base change, and $V$ is an element of
the dense orbit over an algebraic closure precisely when $\mExt^1(V,V)=0$. Since $\cR(Q; d)$ is
irreducible (as an affine space), there is at most one rigid representation $M(d)$ of dimension vector
$d$. Conversely, when $Q$ is a Dynkin quiver, in particular, when $Q$ is of type $\Aa$, there are
only finitely many orbits. Consequently, for any dimension vector $d$ up to isomorphism there is
precisely one representation that has dimension vector $d$ and is rigid. Using this fact, we can
define an equivalence relation on the possible dimension vectors as follows. We say $d$ and $d'$
are {\sl equivalent} provided the indecomposable direct summands of $M(d)$ and $M(d')$ coincide
(up to positive multiplicity). Any module $M(d)$ can have at most $n$ pairwise non-isomorphic indecomposable 
direct summands. Let us assume $M(d)$ has just $n-1$ such indecomposable direct summands (thus
it is almost complete). Then there are at most two complements $M_1$ and $M_2$, meaning that $M(d)
\oplus M_1$ is rigid, $M(d)\oplus M_2$ is rigid, and neither $M_1$ nor $M_2$ is already a direct summand
of $M(d)$. Then it is known that (up to renumbering) there is an exact sequence $M_1\lra M\lra M_2$,
with $M$ consisting of direct summands of $M(d)$. Such a sequence is called {\sl exchange sequence}.
In type $\Aa$, $M$ has at most two indecomposable summands. These sequences play a
crucial role for the recursive construction of all tilting modules. In fact, for any equivalence class of
dimension vectors $d$ there exists a unique $d'$ in this class so that $M(d')$ has no multiple 
indecomposable direct summands. Such a module is called {\sl basic}. The maximal ones among
those modules are the tilting modules that contribute to the volume of the polytope $P(Q)$.
All the other ones correspond to certain faces. We are dealing with the maximal ones, all others
contribute with volume $0$.

The principal aim of this paper is to determine the number of maximal equivalence classes $M(d)$
by using polytopes together with their basic representations. The situation becomes quite elementary
if $Q$ is a quiver of type $\Aa_n$ with its directed orientation. The details can already be found
in \cite{VolHille}. If $Q$ is of type $\Aa_n$ with another orientation, the situation is slightly different,
however the approach in \cite{BaurHille} can be modified as follows. Instead of diagrams with all
connections on the top, we use for arrows from $i+1$ to $i$ connections at the bottom of the
diagram. This even gives a constructive way to compute $M(d)$ for any $d$ and any orientation
of the quiver.

Finally, we define the cones $\sigma$ associated to tilting modules, support tilting modules, and
$2$-support tilting modules, respectively. We start with any module $M$ and its decomposition
into indecomposable direct summands $M=\oplus M(i)^{a(i)}$, where $a(i)$ is the multiplicity of
the indecomposable direct summand $M(i)$ in $M$. For such an $M$ we define $\sigma_M=\mconv
\{\udim M(i), 0 \mid i \in I\}$ to be the convex hull of zero and the dimension vectors of the
indecomposable direct summands. Note that the multiplicities $a(i)>0$ don't play a role in the
definition. Also note that the dimension vectors of the indecomposable direct summands of a tilting
module $T$ form an integral basis of $\Zz^n$, and thus $\sigma_T$ is a simplex which has volume
$1$ by our definition of the volume. If $T$ is a support tilting module, we add the negative standard
basis of the complement of the support of $T$, thus we define $\sigma_T=\mconv\{\udim T(i),-e_j,
0\mid i\in I,j\in J\}$, where $T=\oplus T(i)^{a(i)}$ and $J$ is the complement of the support of $T$.
Finally, for a $2$-support tilting module $T$ we decompose both $T^+$ and $T^-$ as $T^+=
\oplus_{i\in I}T(i)^{a(i)}$, $T^-=\oplus_{j\in J}T(j)^{a(j)}$ and define $\sigma_T=\mconv\{\udim
T(i),-\udim T(j),0\mid i\in I, j\in J\}$.

We illustrate the construction of the simplices and the two series of polytopes by two examples,
namely, quivers of types $\Aa_2$ and $\Aa_3$.
\bigskip

{\sc Example 1.} Let $Q$ be a quiver of type $\Aa_2$. Then there is just one orientation up to a
permutation of the two vertices. The dimension vectors of the indecomposable representations
are $(1,0)$, $(1,1)$, and $(0,1)$. Thus for the roots we get
$$
\Phi^+=\{(1,0),(1,1),(0,1)\}\,,\qquad\Phi^{\mclus}=\{(1,0),(1,1),(0,1),(-1,0),(0,-1)\}\,,\mand
$$
$$
\Phi=\{(1,0),(1,1),(0,1),(-1,0),(-1,-1),(0,-1)\}\,.
$$
The convex hull $C^+(Q)$ of $\Phi^+$ has volume $2$, the convex hull $C^{\mclus}(Q)$ of
$\Phi^{\mclus}$ is a pentagon of volume $5$, and the convex hull $C(Q)$ of $\Phi$ is a hexagon
of volume $6$. The following pairs of roots, together with zero, form a simplex in $P^+(Q)$:
$(1,0),(1,1)$, and $(1,1),(0,1) $. In the polytope $P^{\mclus}(Q)$ we have three additional
simplices defined by the pairs $(0,1),(-1,0)$, $(-1,0),(0,-1)$, and $(0,-1),(1,0)$. The second
of these simplices is replaced in $P(Q)$ by two simplices defined by $(-1,0),(-1,-1)$ and $(-1,-1),
(0,-1)$.
\bigskip

{\sc Example 2.} The case $\Aa_3$ is more complicated, since we have two, essentially different,
orientations in the quiver $Q$. We first describe the parts independent of the orientation, a difference
only occurs for the polytopes of the $P$--series. For the roots we obtain
$$
\Phi^+=\{(1,0,0),(0,1,0),(0,0,1),(1,1,0),(0,1,1),(1,1,1)\}\,,
$$
$$
\Phi^{\mclus}=\Phi^+\cup\{(-1,0,0),(0,-1,0),(0,0,-1)\}\,,\quad\Phi=\Phi^+\cup-\Phi^+\,.
$$
The corresponding convex hulls have volume $5$ for $C^+(Q)$, volume $14$ for $C^{\mclus}(Q)$,
and volume $20$ for $C(Q)$. This can easily be seen from the decomposition of the $P$--series. We
have to chose an orientation for this and describe the triples defining a simplex for the directed orientation
first:
$$
\begin{array}{l}
\{(1,0,0),(1,1,0),(1,1,1)\}\,,\{(0,1,0),(1,1,0),(1,1,1)\}\,,\{(1,0,0),(1,1,1),(0,0,1)\}\,,\\
\{(0,1,1),(0,1,0),(1,1,1)\}\,,\{(1,1,1),(0,1,1),(0,0,1)\}\,.
\end{array}
$$
For the other orientation we get two simplices replaced by two others, however, the union of both
pairs is the same. For this we replace the second and the fourth one by
$$
\{(1,1,0),(0,1,0),(0,1,1)\}\mand\{(1,1,0),(1,1,1),(0,1,1)\}\,.
$$
For the polytope $P^{\mclus}(Q)$ we need to add for any pair of roots, obtained by deleting $(1,1,1)$,
a simplex, and we also need to add for any simple root one simplex. Thus we get for example $\{(1,0,0),
(1,1,0),(0,0,-1)\}$ for the first simplex in $P^+(Q)$, and we get $(1,0,0),(0,-1,0),(0,0,-1)$ for the
simple root $(1,0,0)$. Finally, we add $(-1,0,0),(0,-1,0),(0,0,-1)$ and obtain for the volume $5+5+
3+1=14$. In a similar way we get for the volume of $P(Q)$ the sum $5+5+5+5=20$. Alternatively,
using the computation of the volume with counting the facets and their volume for $C(Q)$, the volume
of $P(Q)$ is $2\times 6+8=20$. This comes from the fact that $C(Q)$ has $6$ squares and $8$ triangles
as facets.

\section{Proofs of Theorems \ref{Thmnumbervol}--\ref{Thmcompare}}

We first prove Theorem \ref{Thmnumbervol}:
It is well known that for any tilting module $T=\oplus T(i)$ the dimension vectors $\udim T(i)$ form a
$\Zz$--basis of $\Zz^n$. Thus $\mvol\sigma_T=1$. Thus we have proven the following lemma.

\begin{Lemma}
For any tilting module, any support tilting module, and any $2$--support tilting module $T$ we have
$\mvol\sigma_T=1$.
\end{Lemma}

It is therefore sufficient to show that $\mvol(\sigma_T\cap\sigma_{T'})=0$ for any two different tilting
modules $T$ and $T'$. This follows from the definition of the volume in \cite{VolHille}.

Crucial for our computation is now Theorem \ref{Thmequal}, which is not true for arbitrary Dynkin quivers.
Note that any non-trivial exchange sequence
$$
0\lra T(i)\lra\oplus T(j)\lra T(i)'\lra 0
$$
has at most two middle terms. Interpreting this as the relation $\udim T(i)+\udim T(i)'=\sum\udim T(j)$,
we see that $P^+(Q)$ is convex, strictly convex at the common facet for one middle term, and flat for
two middle terms. Consequently, $P^+(Q)$ is convex precisely when there are at most two middle terms
for any exchange relation.

Note that such a relation corresponds to two simplices $\sigma_T$ and $\sigma_{T'}$ with a common
facet. If $P^+(Q)$ is convex, then so are $P(Q)^{\mclus}$ and $P(Q)$. Since $P(Q)$ is convex and
has the roots as its vertices, it must coincide with $C(Q)$. The same argument works for $P^+(Q)$
and $P^{\mclus}(Q)$.

The proof of Theorem \ref{Thmcompare} follows from the decomposition of the polytope $P^{\mclus}
(Q)$ with respect to the possible quadrants. Any subset $I$ of the vertices $Q_0$ of $Q$ defines the
quadrant consisting of non-negative entries, whenever the index is not in $I$, and non-positive otherwise.
The volume of $P^{\mclus}(\Aa_n)$ intersected with this quadrant has the same volume as $P^+(\Aa_n)$
intersected with the corresponding face. Thus its volume coincides with the volume of $P^+(\Aa_n|_{Q_0
\setminus I})$. A similar argument holds for $P(\Aa_n)$. Here we need to determine again the volume
in the quadrant defined by $I$. The volume in this case coincides with the product of the volume of
$P^+(\Aa_n|_{Q_0\setminus I})$ and the volume of $P^+(\Aa_n|_I)$.

Theorems \ref{Thmcompute} and \ref{Thmcompute1} can also be proven by using the combinatorial
interpretation in Section 6. Alternatively, one could have used induction over the facets. Since these
proofs need some detailed computations, we defer them to the next section.


\section{Proofs of Theorems \ref{Thmcompute} and \ref{Thmcompute1}}

We first recall the recursion formula $t^+(\Aa_n)=\sum_{i= 0}^{n-1}t^+(\Aa_i)t^+(\Aa_{n-1-i})$.
From this we get one of the standard recursions of the Catalan numbers $C_n$, since $t^+(\Aa_0)
=t^+(\Aa_1)=1$. This is the same recursion as for the number of trees in $B_n$ which will be considered
in the last section. Further details can be found in Stanley's book \cite{Stanley}. We start with the
analogous formula for the polytope $P^+(\Aa_n)$ for the directed orientation. 

\begin{Lemma}
$\mvol P^+(\Aa_n)=\sum\limits_{i=1}^n\mvol P^+(\Aa_n|_{Q_0\setminus\{i\}})$.
\end{Lemma}

This simply uses the fact that there exists precisely one sincere root for $\Aa_n$. Thus, the volume of
$P^+(\Aa_n)$ is just the sum of the volumes of the facets of $P^+(\Aa_n)$ not containing $0$. This
is obvious for the quiver of type $\Aa_n$ with its directed orientation, since every tilting module contains
the projective injective (having the sincere root as dimension vector) as a direct summand. Since $P
(\Aa_n) = C(\Aa_n)$ (independent of the orientation of the arrows) we get the same formula for the
volume. Taking away the projective injective direct summand, we obtain a partial tilting module with
support at $\Aa_n|_{Q_0\setminus\{i\}}$ for precisely one vertex $i$ of the quiver $\Aa_n$. Such a
partial tilting module corresponds to a facet (defined by $d_i=0$). Thus, each tilting module corresponds
to precisely one facet in $P^+(\Aa_n)$ not containing $0$. Moreover, the volume of $\sigma_T$ and
the volume of each of its faces (in particular, of the facet from which we delete the sincere root) is
always $1$. Hence we obtain the formula in Lemma 4.1, and the last formula in Theorem \ref{Thmcompute}
follows directly.

In a next step we compute the volumes of all the polytopes $P^*(\Aa_n)=C^*(\Aa_n)$. This can
be done in several ways. We give a combinatorial approach using certain paths from $(0,0)$ to
$(0,2n)$ later. Firstly, we use the volume of the facets of $C(\Aa_n)$.

\begin{Lemma}
$\mvol C(\Aa_n)=\sum\limits_{i=1}^n{n-1\choose i-1}{n+1\choose i}={2n\choose n}$.
\end{Lemma}

We start by explaining the first equality. Each hyperplane $d_i=1$ contains precisely one facet $F_i$,
and each facet is in the orbit under the symmetric group $S_{n+1}$ (that is, the Weyl goup of the root
system $\Aa_n$) of precisely one such facet. Thus we need to compute the orbit of $F_i$ and the
volume of the facet $F_i$. The volume is just $\left(n-1\atop i-1\right)$ and the orbit has exactly
$\left(n+1\atop i\right)$ elements. This shows the first equality. The second one is a simple recursion
using binomial coefficients:
$$
{2n\choose n}={2n-1\choose n-1}+{2n-1\choose n}={2n-2\choose n-2}+2{2n-2\choose n-1}+
{2n-2\choose n}=\ldots
$$
This proves the first formula in Theorem \ref{Thmcompute}. The formula in Theorem \ref{Thmcompute1}
is just the Catalan recursion, as well as the second formula in Theorem \ref{Thmcompute}. This
finishes the proofs of Theorems \ref{Thmcompute} and \ref{Thmcompute1}.

In the next section we will need another formula which will be used to determine the number of tilting
(or full strongly exceptional) sequences and to relate them to our combinatorial description. Note that the
right hand side is $n+1$ times $n!$.

\begin{Lemma}\label{Loverline}
$(n+1)!=0!n!+1!(n-1)!n+2!(n-2)!\frac{n(n-1)}{2}+\ldots$
\end{Lemma}


\section{Tilting sequences}

If we replace a tilting module by an ordered tuple of modules compatible with non-vanishing homomorphisms,
we even get an easier formula. We define $\overline\cT^+(Q)$ to be the set of tilting sequences, i.e.,
$(T(1),...,T(n))$ satisfying two conditions: $T=\oplus T(i)$ has no self extensions and $\mHom(T(j),T(i))
=0$ for all $j>i$. Note that a tilting sequence is also called a full strongly exceptional sequence of modules.
In a similar way we define support tilting sequences and $2$--support tilting sequences. Then we get the
following formulas for the corresponding numbers, as we will see below.

\begin{Thm}\label{Thmoverline}\label{ThmTiltSequ}
$$
\overline t(\Aa_n)=(n+1)!=0!n!+n1!(n-1)!+\frac{n(n-1)}{2}2!(n-2)!+\ldots\,,
$$
$$
\overline t^+(\Aa_n)=n!\,.
$$
\end{Thm}

The result follows from the interpretation of the map considered in \cite{VolHille}. Define the set $B_n$
as the set of all $3$--regular trees with $n+1$ leaves and one root. There is a natural map $S_n\lra
B_n$ from the symmetric group to the set of all those trees. The number of elements in the preimage
of this map is just the number of tilting sequences (interpreted as a tree with a compatible order on
the inner vertices) defining the same tilting module (interpreted as a tree). Thus we have a bijection
between tilting sequences and elements of the symmetric group.

In order to define the corresponding polytope, we extend the positive roots to $\overline\Phi^+$
consisting of the positive roots for $\Aa_n$ {\sl together with all sums of orthogonal roots}. In case
$Q$ is of type $\Aa_3$ we just add the vector $(1,0,1)$. In terms of dimension vectors we simply
consider vectors with arbitrary entries $0$,$1$ that are non-zero. Using this set as vertices, we can
define $\overline P^+(\Aa_n)$ as the convex hull of zero and $\overline\Phi$. In a similar way, we
define $\overline P(\Aa_n)$ as the convex hull of $\overline\Phi^+$ and $-\overline\Phi^+$, and
$\overline P^{\mclus}(\Aa_n)$ as the convex hull of $\overline\Phi^+$ and the negative simple
roots. To complete the picture, we also need to define a simplex $\overline\sigma_T$ for every
tilting sequence $T$. This can be done as follows. Whenever we have a sequence $(T(i(1)),T(i(2)),
\ldots T(i(r)))$ of indecomposable direct summands of $T$ with $i(1)<i(2)<i(3)<\ldots<i(r)$ and
all components being incomparable (no homomorphisms and no extensions between different members),
then we consider the vertices $\udim T(i(1)),\udim T(i(1))+\udim T(i(2)),\ldots,\udim T(i(1))+\ldots+
\udim T(i(r))$. In this way we get different simplices for different tilting sequences, and the union of all
simplices $\overline\sigma$ of the tilting sequences for a given tilting module $T$ is just the simplex
$\sigma$ of $T$. Thus, the number of tilting sequences is just the volume of $\overline P^+(\Aa_n)$,
the number of support tilting sequences is the volume of $\overline P^{\mclus}(\Aa_n)$, and the
number of $2$--support tilting sequences is the volume of $\overline P(\Aa_n)$. The first and the
last of these numbers have been computed in Theorem \ref{ThmTiltSequ}.

For the corresponding polytope $\overline P^+(\Aa_n)$ we can form $\overline P(\Aa_n)$ just as
the join of $\overline P^+(\Aa_n)$ with $-\overline P^+(\Aa_n)$. This defines a polytope $\overline
P(\Aa_n)$, and its volume is the number of $2$--support tilting sequences. The counting in the above
theorem then computes the volume of $\overline P(\Aa_n)$ from the volume of $\overline P^+(\Aa_n)$
using the formula with the volume of the facets. So we have to compute the facets (in particular, in
all quadrants except the positive and the negative) and their volumes.

The facets in the positive quadrant correspond to elements of the symmetric group, and each facet
contributes with volume $1$. A facet in a quadrant corresponding to a subset $I$ of the vertices
of the quiver corresponds to a facet for $\Aa_n|_I$ and an opposite facet for $\Aa_n|_{Q_0
\setminus I}$. Thus Lemma \ref{Loverline} just computes the volume of $\overline P(\Aa_n)$
from $\overline P ^+(\Aa)$ and Theorem \ref{Thmoverline} is proven. Summarizing this we have
the following result.

\begin{Thm}
$\overline P(\Aa_n)=\overline C(\Aa_n)$ have as volume the number of $2$--support tilting sequences
$(n+1)!$. Moreover, the number of tilting sequences for $\Aa_n$ is $n!$, that is, the volume of
$\overline P^+(\Aa_n)=\overline C^+(\Aa_n)$.
\end{Thm}


\section{Some further comments}\label{sectDyck}

In the final section we present some further explanation for our use of the notion of the volume.
This method is inspired by toric geometry and lattice polytopes. Moreover, we give another interpretation
of the computation of the volume using Stanley's exercise of a combinatorial interpretation of the
Catalan numbers (see Exercise 6.19 in \cite{Stanley}). Very surprisingly, it does not only give an
interpretation of the volume of $P^+(\Aa_n)$ and  $P^{\mclus}(\Aa_n)$ (where the Catalan
numbers occur naturally), but also for the volume of $P(\Aa_n)$, if we modify Dyck paths so that
they correspond to the $2$--support tilting modules.

\subsection{The integral volume}

Note that any cube of the form $[0,1]^n$ has volume one in the euclidean metric and can be
decomposed into $n!$ many simplices, all of volume $1$. This is just a recursive computation.
For $n=$ $1$,$2$ the claim is obvious. Then proceed by induction, and observe that the cube
has $n$ facets containing $0$. (In fact, we could use any vertex instead of $0$). By induction,
the formula holds for the facets, and consequently, for the convex hull of the facet and $(1,
\ldots,1)$ (a pyramid over the facet). Now one checks that the $n$ pyramids over the $n$
facets decompose the cube into $n$ pyramids of volume $(n-1)!$. Thus the cube has volume $n!$.

\subsection{Rooted trees and tilting modules}

We consider the set $B_n$ of rooted $3$--regular trees with one root and $n+1$ leaves (examples
can be found in \cite{VolHille}). If we consider the quiver $\Aa_n$ with its directed orientation, then
we can identify $\cT^+(Q)$ with $B_n$ (see \cite{VolHille}). Thus, we can compute the number of
tilting modules using a standard recursion formula for the number of trees. Take such a tree and
take the unique vertex connected to the root. Decompose the tree $S$ into the two connected
components $S^+$ and $S^-$ obtained from deleting this vertex. Then we get
$$
C_n=\sharp B_n=\sum_{i=0}^{n-1}\sharp B_i\sharp B_{n-1-i}\,,\quad\sharp B_0=\sharp B_1=1\,.
$$
This is one of the standard recursion formulas for the Catalan numbers $C_n$.

\subsection{Dyck paths}

Dyck paths can be used for a combinatorial description of the volume of the polytopes. A {\sl Dyck
path} is a path from $(0,0)$ to $(0,2n)$ using only steps $(1,-1)$ or $(1,1)$ so that the path never
goes below the $x$--axes (meaning that the first coordinate of a point is non-negative). We denote
the set of Dyck paths by $D^+_n$. For given $n$ the number of Dyck paths coincides with the
number of possible bracketings of an expression with $n$ inputs. Moreover, this can be identified
with the elements in $B_n$ and with the vertices of the associahedron. For the combinatorics we
refer to the famous exercise in Stanley's book \cite{Stanley}, where we use only 5 interpretations
of the $66$ (in fact, even more can be found on Stanley's homepage). If we consider arbitrary paths
from $(0,0)$ to $(0,2n)$ with steps $(1,1)$ or $(1,-1)$ (without the condition to be above the $x$--axes)
we obtain a set $D_n$ that has $2n\choose n$ elements. An interpretation of $2$-support tilting sequences
is obtained as follows. Whenever the path stays above the $x$--axes, we identify the corresponding
Dyck path with its tree, and thus with a direct summand $T^+$ of $T$. Whenever the path stays below
the $x$--axes, we identify the corresponding path with $T^-$. This bijection identifies paths in $D_n$
with $2$--support tilting modules for $\Aa_n$. Consequently, we have computed the volume of
$P(\Aa_n)$ as the number of elements in $D_n$ which is $2n \choose n$.


\bibliographystyle{amsalpha}

\end{document}